\documentclass[a4paper,11pt]{article}
\usepackage{amsfonts,amssymb,amsmath,amsthm,latexsym,amscd,graphicx,xy,psfrag}
\usepackage[dvipsnames]{xcolor}
\usepackage{epic, eepic,booktabs}
\usepackage[colorlinks=true]{hyperref}
\hypersetup{urlcolor=blue, citecolor=red}

\def\<#1,#2>{\langle #1,#2 \rangle}

 \def\bothID{\rlap{\hbox to.97\wd0{\hss\vrule height.06\ht0 width.82\wd0}}
 \copy0\rlap{\kern-.36\wd0\vrule height1.05\ht0 width.05\ht0}\kern.14\wd0}

 \DeclareMathOperator{\spec}{spec}
 
\DeclareMathOperator{\vol}{vol}

 \DeclareMathOperator{\Spin}{Spin}
 \DeclareMathOperator{\SO}{SO}

\DeclareMathOperator{\re}{Re}

\newcommand{\Dt}{D^{N_t}}
\newcommand{\Duu}[1]{D^{N_{#1}}}

\newcommand{\nat}{\nabla^{N_t}}

\newcommand{\nuu}{\nu}
\newcommand{\pt}{\Pi_{t_0}^{t}}

\newcommand{\h}{g|_{\partial N}}

\begin{document}

\title{Riemannian Metrics and Harmonic Sections of Spinor Bundles}
\author{Simone Farinelli\\
        Aum\"ulistrasse 20\\ CH-8906 Bonstetten\\Email:
        simone@coredynamics.ch}
\maketitle
\begin{abstract}
We study the clustering of the lower eigenvalues of the
Dirac operator on a spinor bundle when the metric structure
is varied. We then apply the general results to show that any closed spin
manifold of dimension $m\ge4$ has a Riemannian metric admitting non-trivial harmonic spinors.

\end{abstract}

\newtheorem{theorem}{Theorem}[section]
\newtheorem{proposition}[theorem]{Proposition}
\newtheorem{lemma}[theorem]{Lemma}
\newtheorem{corollary}[theorem]{Corollary}
\theoremstyle{definition}
\newtheorem{ex}{Example}[section]
\newtheorem{rem}{Remark}[section]
\newtheorem*{nota}{Notation}
\newtheorem{defi}{Definition}
\newtheorem*{conjecture}{Conjecture}
\renewcommand{\P}{P}

\section{Introduction}
On a closed connected $m$-dimensional Riemannian spin manifold $(M, g, \P)$ with metric
$g$ and spin structure $\P=\P_{\mathrm{Spin}}(M)$ we consider the Dirac operator
$D$, a first order differential operator which acts on sections of the spinor bundle $\Sigma$. It is known that the dimension $h(M,g,\P)=\dim\ker(D)$ of the space $\ker(D)$ of harmonic spinors is not a topological invariant. For instance, the
one parameter family of Berger metrics on $S^3$ admits non-trivial harmonic spinors only for non-generic values of the parameter, cf. \cite{Hi74}.

There are, however, topological bounds for the metric invariant $h(M,g,P)$.
In even dimension $\Sigma$ splits into the
direct sum $\Sigma=\Sigma^{+} \bigoplus \Sigma^{-}$ of the bundles of positive and negative Weyl spinors and $D$ is of the form
\begin{equation}D=\left[\begin{matrix} 0 & D^- \\ D^+ & 0 \end{matrix}\right].
\end{equation}
We write $h(M,g,\P)=h^{+}(M,g,\P)+h^{-}(M,g,\P)$ where $h^{\pm}(M,g,\P):=\ker(D^{\pm})$.
If $m=4k$, by the Atiyah-Singer index Theorem
 \begin{equation}h^{+}(M,g,\P)-h^{-}(M,g,\P)=\widehat{A}(M),\end{equation}
where $\widehat{A}(M)$ is the $\widehat A$-genus of $M$, cf. \cite{AS68}, therefore $h(M,g,\P)\ge|\widehat{A}(M)|$ for \textit{any} Riemannian metric and spin structure.
There is a similar index thereom if $m\equiv 1,2 \mod 8$, cf. \cite{Mi65, AS71}:
\begin{equation}
   \begin{split}
           m=8k+1\Rightarrow  h(M,g,\P) &\equiv
           \alpha(M)\mod 2,\\
           m=8k+2 \Rightarrow  \frac{h(M,g,\P)}{2} &\equiv
         \alpha(M) \mod 2,
   \end{split}
 \end{equation}
where $\alpha(M)\in\mathbf Z_2$ is Milnor's $\alpha$-genus, a topological invariant which depends on the choice of the spin structure on $M$. In particular $h(M,g,P)\neq 0$ if $\alpha(M)\neq 0$.

In more recent years, there have been extensive studies on $D$-minimal Riemannian metrics, that is metrics for which the lower bound given by the index theorem is attained:

\begin{center}
\begin{tabular}{|l|l|l|}\hline
  $\dim(M)$ & Genus & $h(M,g,P)$\ for $D$-minimal metrics \\ \hline\hline
  $m=4k$ &  $\begin{gathered} \\ \widehat A(M)\ge 0\end{gathered}$ &  $h^{+}(M,g,P)=\widehat{A}(M)$\\
  &   &   $h^{-}(M,g,P)=0$\\\cline{2-3}
  &   $\begin{gathered} \\ \widehat A(M)< 0\end{gathered}$ &  $h^{+}(M,g,P)=0$\\
  & &  $h^{-}(M,g,P)=-\widehat{A}(M)$\\\hline
  $m=8k+1$ & $\alpha(M)=0$  &  $h(M,g,P)=0$\\\cline{2-3}
  & $\alpha(M)=1$  &  $h(M,g,P)=1$\\ \hline
  $m=8k+2$ & $\alpha(M)=0$ &  $h(M,g,P)=0$\\\cline{2-3}
   & $\alpha(M)= 1$&  $h^{+}(M,g,P) =h^{-}(M,g,P)=1$\\\hline
  $m\equiv 3,5,6,7\mod 8$ & $\widehat A(M)=0^{\phantom{C^{C^{C}}}}$ &  $h(M,g,P)=0$\\ \hline
\end{tabular}
\end{center}

In \cite{Ma97} it is proved using a variational approach that any generic metric on a closed connected spin manifold of dimension $\leq 4$ is D-minimal. In
\cite{BD02} the same result is proved for simply connected manifolds of dimension at least $5$, using surgery and spin bordism arguments.
These results were later established in \cite{ADH09} for any dimension  and without the simply connectedness hypothesis; the essential step in the proof is the construction of a single $D$-minimal metric on any given closed connected spin manifold.

On the other hand, there is the following long-standing conjecture, cf. f.i. \cite{BD02}.
\begin{conjecture}
Let $(M,P)$ be a closed connected spin manifold of dimension $\ge3$. Then, there exists a Riemannian metric $g$ on $(M,P)$ with non-trivial harmonic spinors, i.e. $h(M,g,P)>0$.
\end{conjecture}
The conjecture does not hold in dimension $2$, because
for any eigenvalue $\lambda$ of the Dirac operator on a closed surface $\Sigma$ of genus zero,  one has
\begin{equation}\lambda^2\ge\frac{4\pi}{\vol(\Sigma)},
\end{equation}
in particular, there are no harmonic spinors (for any metric and spin structure), cf. \cite{Ba92}. If the genus is $\le 2$ then all metrics are $D$-minimal and the spectrum of $D$ can be explicitly computed when the genus is $1$, cf. \cite{BS92}. For results on spin structures for all hyper-elliptic surfaces $\Sigma$ of genus $\ge 2$ we refer the reader to \cite{BS92}.

The conjecture has been shown to be true for $m\equiv 0,1,7\mod 8$ by Hitchin, using the Atiyah-Singer index theorem and the theory of exotic spheres in \cite{Hi74},
and for $m\equiv 3,7\mod 8$ by B\"ar using an analytic approach in \cite{Ba96}.
Both these results follow also from Dahl's theorem on invertible Dirac
operators, cf. \cite{Da08}.

For the special case of $m$-spheres, explicit metrics with non-trivial harmonic spinors
have been constructed when $m\equiv 0\mod 4$ in \cite{Se00} and in the case of Berger spheres when $m\equiv 3\mod 4$ in \cite{Ba96}.

In particular Dahl proves in \cite{Da05} that for any closed connected spin manifold of dimension $m\ge 3$ and any {\it non-zero} $\lambda\in\mathbf R$ it is possible to construct a Riemannian metric for which the Dirac operator has $\lambda$ has eigenvalue. The main contribution of this paper is to extend Dahl's result to $\lambda=0$ and establish that the conjecture holds in all dimensions $m\ge 4$:
\begin{theorem}\label{thm11}
Let $(M,P)$ be a closed connected spin manifold of dimension $m\ge4$. Then, there exists a Riemannian metric $g$ on $(M,P)$ admitting non-trivial harmonic spinors.
\end{theorem}
This paper is structured as follows. Section 2 recalls the basics of
the theory of Dirac bundles. In Section 3 we establish several results on the spectrum of Dirac operators, preliminary to the proof of the main theorem in Section 4. Section 5 concludes
with some observations.

\section{Preliminaries}
\subsection{Dirac bundles}
The purpose of this section is to recall the basic definitions and set the notations
concerning the theory of Dirac operators. General references are \cite{LM89}, \cite{BW93} and \cite{BGV96}.
\begin{defi}\label{DiracBundle}
Let $\Sigma\to M$ be a complex vector bundle on a Riemannian manifold endowed with a Hermitian structure $\<\cdot,\cdot>$ and a Hermitian connection $\nabla$. Then $\Sigma$ is a \textbf{Dirac bundle} if it is a bundle of Clifford modules such that:
\begin{enumerate}
    \item[(i)] Clifford multiplication $X\cdot$ by vectors is skew-adjoint w.r.t. $\<\cdot,\cdot>$;
   \item[(ii)] Clifford multiplication is compatible with $\nabla$, i.e.
\begin{equation}
\begin{split}
\nabla(X\cdot\varphi)=\nabla X\cdot\varphi+X\cdot\nabla\varphi
\end{split}
\end{equation}
for all $\varphi\in\Gamma(\Sigma)$.
\end{enumerate}
\end{defi}
Given a Dirac bundle $\Sigma$, the associated \textbf{Dirac operator} $D:\Gamma(\Sigma)\to \Gamma(\Sigma)$ is the first-order differential operator defined by
the composition
\begin{equation}
\begin{CD}
   {\Gamma(\Sigma)} @>{\nabla}>> {\Gamma(T^*M\otimes \Sigma)}@>{}>>\Gamma(\Sigma)\;,
   \end{CD}\end{equation}
of $\nabla$ with Clifford multiplication.
The square $D^2$ of the Dirac operator is the \textbf{Dirac Laplacian}.
\vskip0.2cm\par
A particularly important example of Dirac bundle is given by the spinor bundle.
We recall that an oriented Riemannian manifold $(M,g)$ of dimension $m\geq 3$ is a \textbf{spin manifold} if it admits a spin bundle, i.e.~a double cover $\pi:P\to SO(M)$ of the oriented orthonormal frame bundle $SO(M)$ of $M$ with structure group $\Spin(m)$ inducing the canonical covering $\Theta:\Spin(m)\to \SO(m)$ on each fiber.
The vector bundle
      \begin{equation}
      \Sigma=P\times_{\Spin(m)}\mathbf{C}^l\;,\qquad
      l=2^{[\frac{m}{2}]}\;,\end{equation}
associated with the spin representation $\mathbf C^l$ of $\Spin(m)$ is the \textbf{spinor bundle}. It is a Dirac bundle with respect to the Levi-Civita connection $\nabla$ on $M$.

The restriction of $\Sigma$ to an oriented hypersurface $N\subset M$ is a Dirac bundle, as we now briefly recall. Let $\nu$ be the unit normal vector field and $B:TN\to TN$ the form operator with respect to $\nu$, $B(X)=-\nabla_X\nu$. We have the classical Gauss formula
$$
\nabla_XY=\nabla^N_XY+g(B(X),Y)\nu
$$
for all $X,Y\in\mathfrak X(N)$, where $\nabla^N$ is the Levi-Civita connection on $(N,g|_N)$.
\begin{proposition}[see e.g. \cite{Ba96}]\label{DiracBoundary}
Let $\Sigma$ be the spinor bundle on a
Riemannian manifold $M$ and $N\subset M$ an oriented
hypersurface. Then $\Sigma|_N$ has a natural structure of Dirac bundle on $N$ where
\begin{equation*}
\begin{aligned}
X\bullet\sigma&=X\cdot\nu\cdot\sigma\;,\\
\nabla_X\sigma&=\nabla^N_X\sigma+\frac{1}{2}B(X)\cdot\nu\cdot\sigma
\end{aligned}
\end{equation*}
for all $X\in TN$ and $\sigma\in\Gamma(\Sigma|_N)$.
\end{proposition}
\subsection{Spectral Resolutions}
\label{Spectrum}
We review here the spectral properties of the  Dirac operator and Laplacian over a compact Riemannian manifold $M$. Consider first the case $\partial M=\varnothing$, see e.g.~\cite{BW93,Gi95}. The operators $D, D^2:\Gamma(\Sigma)\to\Gamma(\Sigma)$ are elliptic and formally self-adjoint; taking as domains the closure of $\Gamma(\Sigma)$ in the Sobolev $H^1$- and $H^2$- topology leads to two self-adjoint operators: 

\begin{proposition}\label{PQwithoutBoundary}
For compact $M$ the operator $D$ admits a discrete spectral resolution, i.e., there exists a sequence $(\varphi_j,\lambda_j)_{j\ge0}$ such that
$(\varphi_j)_{j \in \mathbf{N}}$ is an orthonormal basis of $L^2(M,\Sigma)$ consisting of smooth sections and, for all $j\ge0$,
\begin{equation}
 D \varphi_j=\lambda_j \varphi _j,
\end{equation}
with $\lambda_j\in\mathbf{R}$. The absolute values of the sequence $(\lambda_j)_{j\ge0}$ can be ordered in an increasing  diverging sequence.
An analogous result holds for $D^2$, in this case $\lambda_j\geq 0$ for all $j\ge0$.
\end{proposition}

The existence of a discrete spectral resolution for Dirac and Laplacian operators on manifolds with boundary is more involved. In particular, while for the Dirac Laplacian
is always possible to find local elliptic boundary conditions
allowing for a discrete spectral resolution, this is not always the case
for the Dirac operator \cite[\S 1.11.6]{Gi95}, see also \cite{Gr96,Ho85}.

The Dirac Laplacian on a Dirac bundle $\Sigma\to M$ over a compact Riemannian
manifold $M$ with boundary is a formally self-adjoint operator on the smooth sections satisfying the Dirichlet boundary condition
$B_D\varphi:=\varphi|_{\partial M}=0$ or the Neumann
boundary condition $B_N\varphi=\nabla_{\nu}\varphi=0$, where $\nu$ is the inward normal vector field on $\partial M$. Taking the closure in the Sobolev
$H^2$-topology, leads to a self-adjoint operator.

The following analogue of Proposition \ref{PQwithoutBoundary} is a special case of classical elliptic boundary theory developed by Seeley in \cite{Sl66, Sl69} and Greiner in \cite{Gr70, Gr71}. Therein
one fixes $B=B_D$ or $B=B_N$.

\begin{proposition}\label{PND}
There exists an orthonormal basis $(\varphi_j)_{j\ge0}$
of $L^2(M,\Sigma)$ consisting of smooth sections and such that, for all $j\ge0$,
\begin{itemize}
\item[(i)] $D^2\varphi_j=\lambda_j \varphi_j$,
\item[(ii)] $B\varphi_j=0$,
\end{itemize}
with $\lambda_j\in\mathbf R$. The absolute values of the eigenvalues $(\lambda_j)_{j\ge0}$ can be ordered in an increasing diverging sequence, with only a finite number of the $\lambda_j$ negative.
The Dirichlet eigenvalues are all strictly positive. The
Neumann eigenvalues are all but for a finite number strictly
positive.
\end{proposition}
\noindent For an example of negative Neumann eigenvalues see \cite{Fa98} Chapter 5. For an extensive overview of Dirac and Dirac Laplacian spectra on manifolds without and with boundary and their elliptic boundary conditions see \cite{Fa23} subsections 3.1 and 3.2.
\subsection{Spectral Estimates}
We review the Courant-Hilbert Theorem,
originally formulated for the scalar Laplacian on $\mathbf{R}^m$ \cite{CH93, Cha84} and  extended to the Dirac Laplacian in \cite{Ba91, Fa98}.

\begin{proposition}[\textbf{Spectral Upper Bounds}]\label{thm27bis}
Let $\Sigma\to M$ be a Dirac bundle over a compact, oriented, Riemannian manifold $M$ with
(smooth) boundary $\partial M$. For some integer $N\geq 1$, consider the decomposition
\begin{equation}M={\bigcup}^N_{k=1}M_k\end{equation}
of $M$ into $0$-codimensional submanifolds $M_k$ with (smooth) boundaries $\partial M_k$ and pairwise disjoint interiors. We require that if a boundary $\partial M_k$
intersects $\partial M$ then it agrees with the corresponding connected component of $\partial M$. Let $(\lambda_j)$ be the Dirichlet spectrum of $D^2$ on $M$, likewise $(\mu_j^k)$ for $M_k$.

Choose any $N'\le N$ and set
\begin{equation}
(\mu_j)=\bigcup_{k=1}^{N'}(\mu_j^k)\,,
\end{equation}
where the sequences are in non-decreasing order and the eigenvalues are counted
according to their multiplicities. Then
\begin{equation} \lambda_j\le\mu_j,\end{equation}
for all $j\ge0$.
\end{proposition}

\par There exists  an analogous result for the Neumann boundary condition \cite{Ba91, Fa98}, but we will not need it in the following.
\section{Warped products and Variations of metrics}
\subsection{Warped Products}
In this section we compute the eigenvalues of the Dirac Laplacian under the
Dirichlet boundary condition on a cylindrical manifold $M:=[0,1]\times N$, where $(N,g_N)$ is a closed and oriented Riemannian manifold of dimension $\dim N=m-1$. We equip $M$ with the warped product metric
\begin{equation}
g:=dt^2+\rho^2(t)g_N,
\end{equation}
where $\rho$ is a smooth function, and assume that $M$ is a spin manifold with spinor bundle $\Sigma\to M$.

By Proposition \ref{DiracBoundary} we have an induced Dirac bundle structure on each hypersurface $(N_t:=\{t\}\times N,\rho^2(t) g_N)$. We note that $\nuu=\frac{\partial}{\partial t}$ is the unit normal vector to any $N_t$ and denote by $\nabla$ and $\nat$
the Levi-Civita connections on $M$ and $N_t$. We now compare the Dirac and Laplacian operators $D$, $D^2$ on $M$ with the analogous operators $\Dt$, $(\Dt) ^2$ on $N_t$.

By \cite{Ba96} one has
\begin{equation}
\label{dirsubm}
D\varphi=\nuu \cdot \Dt \varphi-\frac{m-1}{2}H\, \nuu\cdot\varphi+\nuu \cdot\nabla_{\nuu}\varphi
\end{equation}
for $\varphi\in\Gamma(\Sigma)$. Here and in the following $H= -\rho^\prime/\rho$ is the mean curvature of $N_t$. The following lemma deals with the Laplacians.
\begin{lemma}\label{lemma54}
Let $\varphi\in\Gamma(\Sigma)$ then
\begin{equation}
\begin{split}
D^2\varphi &= (\Dt)^2 \varphi+[\Dt,\nabla_{\nuu}]\varphi
+\left(\frac{m-1}{2}H^{\prime}-\left(\frac{m-1}{2}\right)^2 H^2\right)\varphi\\
&\;\;\;+(m-1)H\,\nabla_{\nuu}\varphi-\nabla_{\nuu}^2\varphi\;.
\end{split}
\end{equation}
\end{lemma}
\begin{proof}
It follows from a straightforward computation making use of (\ref{dirsubm}), the fact that $\nabla_{\nuu}\nuu=0$ and the relations
\begin{equation}
\Dt (\nuu \cdot \varphi)= -\nuu\cdot \Dt\varphi,\;\;
-\nuu\cdot\nabla_{\nuu}(H\nuu\cdot\varphi)
= H^\prime\varphi+H\nabla_{\nuu}\varphi\;,
\end{equation}
which hold for all $\varphi\in\Gamma(\Sigma)$.
\end{proof}

For any $n\in N$ let $\pt $ be parallel transport in $\Sigma$ with respect to $\nabla$ along the curve $u \mapsto (u,n)$ from $(t_0,n)$ to $(t,n)$, where $t_0,t\in[0,1]$. For the sake of brevity, we omit the proof of the following useful lemma.
\begin{lemma}\label{propA}
We have
\begin{equation}\label{paralleltransport}
\Dt \pt\sigma
=\frac{\rho(t_0)}{\rho(t)}\pt\Duu{t_0}\sigma,
\end{equation}
for all $\sigma\in\Gamma(\Sigma|_{N_{t_0}})$.
\end{lemma}
Let $\{\sigma_j\}_{j\ge 0}$ be an $L^2$-orthonormal eigenbasis of the Dirac operator $\Duu{t_0}$ at a fixed slice $N_{t_0}$ with eigenvalues $\{ \mu_j\}$. We note that $\nuu\cdot\sigma_j$ is an eigenvector
with eigenvalue $-\mu_j$, as Clifford multiplication by $\nuu$ anticommutes with $\Duu{t_0}$. We may hence assume $\sigma_{j+1}=\nuu\cdot\sigma_j$, $\mu_{j+1}=-\mu_j$ for, say, all $j$ even.

Parallel transport along $t$-lines yields sections
$$\varphi_j=\pt\sigma_j$$ of the spinor bundle $\Sigma$ on $M$.
By Proposition \ref{propA}, the collection $\{\varphi_j^t\}_{j\ge 0}$ of their restrictions $\varphi_j^t=\varphi_j|_{N_t}$ to $N_t$
 is an $L^2$-orthonormal eigenbasis of $\Dt$ with the eigenvalues
 $$
 \mu_j(t)=\frac{\rho(t_0)}{\rho(t)}\mu_j\;.
 $$
Clearly $\varphi_j^{t_0}=\sigma_j$ and $\varphi_{j+1}^t=\nuu\cdot\varphi_j^t$ at all $t$, since $\nu\cdot$ is parallel along $t$-lines.

Let now $\varphi$ be a smooth section over $M$. Restriction to $N_t$ yields a smooth
section over $N_t$ which can be expressed in the basis $\{\varphi_j^t\}_{j\ge 0}$, therefore
\begin{equation}\label{dec}
\varphi
=\sum_{j\ge 0}a_j\varphi_j\;,
\end{equation}
for some functions $a_j=a_j(t)$.
We use (\ref{dec}) to rewrite the eigenvalue
equation and boundary condition for the Dirac Laplacian on $M$ as follows.
\begin{proposition}\label{prop55}
The equation $D^2 \varphi = \lambda \varphi$ subject to the Dirichlet boundary condition
$\varphi|_{\partial M}=0$ is equivalent to the system
\begin{equation}\label{rep}
\begin{split}
-&a_j^{\prime\prime}+(m-1)Ha_j^{\prime}+\left(\mu_j^2-\mu_j^{\prime}+\frac{m-1}{2}H^{\prime}-\frac{(m-1)^2}{4}H^2-\lambda\right)a_j=0,\\
&a_j(0)=a_j(1)=0
\end{split}
\end{equation}
for all $j\ge 0$. If $N$ admits a non-trivial harmonic spinor, then the equation has a non-trivial harmonic spinor for any eigenvalue of the form $\lambda=\pi^2 n^2$, where $n$ is some positive integer.
\end{proposition}
\begin{proof}
Inserting the decomposition (\ref{dec}) into the equation
$(D^2-\lambda)\varphi=0$ and using Lemma \ref{lemma54}
and $\nabla _{\nuu} \varphi_j=0$,
we obtain
\begin{equation}
\begin{split}
\sum_{j\ge0}&\left[-a_j^{\prime\prime}+(m-1)Ha_j^{\prime}+\right.\\
&+\left.\left(\mu_j^2-\mu_j^{\prime}+\frac{m-1}{2}H^{\prime}-\frac{(m-1)^2}{4}H^2-\lambda\right)a_j\right]\varphi_j=0
\end{split}\end{equation}
and the first equation of (\ref{rep}) follows immediately. The second equation is just the Dirichlet boundary condition.

Now, under the substitution
\begin{equation}
\bar{a}_j:=\exp\left [-\frac{m-1}{2}\int_0^{t}\!\!\!\! H\right]a_j,
\end{equation}
the system (\ref{rep}) becomes
\begin{equation}\label{sys}
\begin{split}
-&\bar{a}_j^{\prime\prime}+\left(\mu_j^2-\mu_j^{\prime}-\lambda\right)\bar{a}_j=0,\\
&\bar{a}_j(0)=\bar{a}_j(1)=0
\end{split}
\end{equation}
for all $j\ge 0$. If $N$ admits a non-trivial harmonic spinor then we have $\mu_k=0$ for some $k$ and hence $\mu_k(t)=0$ identically. We set $\bar{a}_{j}=0$ for all $j\neq k$ and reduce system (\ref{sys}) to
\begin{equation}\label{sys2}
\begin{split}
-&\bar{a}_j^{\prime\prime}-\lambda\bar{a}_j=0,\\
&\bar{a}_j(0)=\bar{a}_j(1)=0,
\end{split}
\end{equation}
which can be easily seen to have non-trivial solutions with $\lambda$ as required.
\end{proof}

\subsection{Continuity of the Eigenvalues}
We recall that the uniform $C^k$-topology for sections of a vector bundle $E\to M$ over a compact manifold $M$ is usually introduced by means of a good presentation of $E$ and a partition of unity, which allow to globalize the corresponding topology for (locally defined) functions, see e.g. \cite{LM89}.

For our purposes, however, it is more convenient to consider the following equivalent definition. If $(M,g)$ is a compact Riemannian manifold and $E$ an Hermitian vector bundle with a compatible connection, we set
\begin{equation}
\label{eq:uniform}
\|\sigma\|^2_{C^k}=\sup_{x\in M}(\sum_{j=0}^k|\underbrace{\nabla\nabla\ldots\nabla}_{j-\text{times}}\sigma|^2)
\end{equation}
for all $\sigma\in\Gamma(E)$, where $\nabla$ is the combination of the connection on $E$ and the Levi-Civita connection of $(M, g)$ and $\left|\cdot\right|$ is the pointwise norm on the fibres of $TM$ and $E$.

The norm \eqref{eq:uniform} {\it does} depend on the choice of metrics and connections but different norms are all equivalent and induce the uniform $C^k$-topology. For $k=0$, this is the
$L^\infty$-topology and we write $\|\sigma\|_{L^\infty}=\|\sigma\|_{C^0}$ for the corresponding norm. Crucially for our aims, this norm still depends on the choice of metrics.\par
Nowaczyk proves in \cite{No13} the following result.
\begin{theorem}\label{ev_cont} Let $M$ be a compact spin manifold and $\mathcal{R}(M)$ be the normed space of the Riemannian metrics over $M$ with the norm inducing the uniform $C^1$-convergence over $M$. For any choice of the Riemannian metric $g\in\mathcal{R}(M)$ let $D^g$ be the Dirac operator over $M$. The spectrum of $D^g$ is a pure point spectrum which can be ordered as
\begin{equation}
\spec(D^g)=(\lambda_j(g))_{j\in\mathbf{Z}},
\end{equation}
where the eigenvalues are repeated according to their multiplicities, and for all $j\in\mathbf{Z}$ and all $g\in\mathcal{R}(M)$
\begin{equation}
\lambda_j(g)\le\lambda_{j+1}(g).
\end{equation}
All the eigenvalues $\lambda_j:\mathcal{R}(M)\rightarrow \mathbf{R}$ are continuous with respect to the uniform $C^1$-topology.
\end{theorem}

\begin{rem}
There are different possibilities of ordering the spectrum:
\begin{itemize}
\item If we fix the order in such a way that $\lambda_0(g)$ is the smallest non-negative eigenvalue near to $0$, then in general the $\lambda_j$'s will not be continuous functions of the Riemannian metric with respect to uniform $C^1$-convergence. Nowaczyk constructs an explicit counter-example in the proof of his Main Theorem 3 in \cite{No13}.

\item Nowaczyk introduces an enumeration for the eigenvalues such that the $\lambda_j$'s are continuous functions of the Riemannian metric with respect to uniform $C^1$-convergence, see the proof of his Main Theorem 2 in \cite{No13}.

\item If we fix the ordering such that the eigenvalues are repeated according to their multiplicity, and for all $j\in\mathbf{Z}$ and all $g\in\mathcal{R}(M)$
    \begin{equation}
      0\le\lambda_j^2(g)\le\lambda_{j+1}^2(g),
    \end{equation}
     the $\lambda_j^2$'s and the $\lambda_j$'s are continuous functions of the Riemannian metric with respect to uniform $C^{\infty}$-convergence, as proved by Canzani in the proof Theorem 2.7 in \cite{Ca14}.
\end{itemize}
\end{rem}

\subsection{Differentiability of the Eigenvalues and  First Variation }
We recall that there exists a geometric process to compare spinor fields for two different Riemannian metrics \cite{BG92}. By a result of Rellich, for an {\it analytic} variation of metrics with associated one-parameter family $(\Sigma^t, \langle \cdot,\cdot\rangle^t,\nabla^t)$ of spinor bundles,
there is an analytic discrete spectral resolution $(\varphi_j^t,\lambda_j^ t)_{j \ge 0}$ of the corresponding family of Dirac operators. As usual, the absolute values of the sequence $(\lambda_j^t)_{j\ge 0}$ can be ordered in an increasing diverging sequence.

Explicit formulas for the first derivative of the analytic branches of the eigenvalues are provided by the following.
\begin{theorem}\label{BGthm}
If $g^t=g+tk$ is a linear variation of the Riemannian metric $g$ on
a compact $m$-dimensional spin manifold $M$, where $k$ is some symmetric $2$-tensor,
then any eigenvalue $\lambda_j(g^t)$ of the Dirac operator $D^{g^t}$ of multiplicity $\mathfrak{m}$ can be written as
\begin{equation}\label{ev_an}
\lambda_j(g^t)=\left\{
                 \begin{array}{ll}
                   \mu_j^p(t), & (t\ge0) \\
                    &\\
                   \mu_j^q(t), & (t\le0)
                 \end{array}
               \right.
\end{equation}
for appropriate $p,q\in\{1,2,\dots,\mathfrak{m}\}$, where $\{\mu_j^1,\dots,\mu_j^\mathfrak{m}\}$ are real analytic real function of $t$ in an open real neighbourhood of $0\in\mathbf{R}$. Moreover, the directional derivatives of the $\mu_j^p$'s at $g$ in direction $k$ read
\begin{equation}
\label{extendedBG}
\begin{aligned}
\left.\frac{d}{dt}\right|_{t=0}\mu_j^p(g^t)&=-\frac{1}{2}\int_M\,g\left( Q_{\varphi_j^0}, k\right)\mathrm{vol}_{g}\\
&=-\frac{1}{2}\int_M\,\left\langle\sum_{i=1}^m e_j\cdot\nabla_{K_g(e_j)}\varphi^0_j,\varphi^0_j\right\rangle \mathrm{vol}_{g}\;,
\end{aligned}
\end{equation}
where $(e_j)$ is a local $g$-orthonormal frame, $Q_{\varphi}$ is the symmetric $2$-tensor
\begin{equation}
\label{eq:Qvar}
Q_{\varphi}(X,Y)=\frac{1}{2}\re\left\langle X\cdot\nabla_Y\varphi+Y\cdot\nabla_X\varphi,\varphi\right\rangle\end{equation}
defined for any $\varphi\in\Gamma(\Sigma^0)$, and  $K_g$ is the endomorphism of $TM$ defined by $k(X,Y)=g(K_g(X),Y)$.
\end{theorem}

\begin{proof}
 By Theorem \ref{Rellich} (Theorem VII.3.9 (Rellich) in \cite[pp.~392--393]{Ka80}), for any symmetric $2$-tensor $k$, the eigenvalue $\lambda_j^t=\lambda_j(g+tk)$ has real analytic branches $\mu_j^1(t),\dots,\mu_j^\mathfrak{m}(t)$, and by adapting Theorem 2.1 in \cite{EI08} from the Laplace-Beltranmi operator to the Dirac operator we obtain equation (\ref{ev_an}). Thereby, we have utilized the continuity of $\lambda_j^t$ at $t=0$. See also Subsection 2.1 in \cite{KMP23}. Formula (\ref{extendedBG}) is proved in \cite[pp.~593--595]{BG92}.\\
\end{proof}

\begin{corollary}\label{ev_diff} If $g^t=g+tk$ is a linear variation of the Riemannian metric $g$ on
a compact $m$-dimensional spin manifold $M$, where $k$ is some symmetric $2$-tensor,
then any eigenvalue $\lambda_j(g^t)$ of the Dirac operator $D^{g^t}$ of multiplicity $\mathfrak{m}$ can be written as
\begin{equation}\label{ev_an}
\lambda_j(g^t)=\left\{
                 \begin{array}{ll}
                   \tilde{\mu}_j^p(g^t), & (t\in[0,+\varepsilon[) \\
                    &\\
                   \tilde{\mu}_j^q(g^t), & (t\in]-\varepsilon, 0])
                 \end{array}
               \right.
\end{equation}
for appropriate $p,q\in\{1,2,\dots,\mathfrak{m}\}$, where $\{\tilde{\mu}_j^1,\dots,\tilde{\mu}_j^\mathfrak{m}\}:\mathcal{R}(M)\rightarrow\mathbf{R}$ are infinitely Fr\'{e}chet differentiable functions for an $\varepsilon>0$. In particular, if $\lambda_j(g)$ is a simple eigenvalue of $D^g$, then it is infinitely Fr\'{e}chet differentiable.
\end{corollary}

\begin{proof}[Proof of Corollary \ref{ev_diff}]\text{}\\
Let us consider the proof of Theorem \ref{BGthm}.
\begin{itemize}
\item By Theorem \ref{Rellich} (Theorem VII.3.9 (Rellich) in \cite[pp.~392--393]{Ka80}), for any symmetric $2$-tensor $k$ the eigenvalue $\lambda_j^t=\lambda_j(g+tk)$ has real analytic branches $\{\mu_j^1,\dots,\mu_j^\mathfrak{m}\}$, which are functions of $t$ on an open neighbourhood $]-\varepsilon, +\varepsilon[$ of $0\in\mathbf{R}$. In particular, for any $g\in\mathcal{R}(M)$ and any direction $k\in T\mathcal{R}(M)$ all the directional derivatives of any order of the branches $\mu_j^p$, for any $p \in\{1,2,\dots,\mathfrak{m}\}$, at $g$ in direction $k$ are well defined for all $j\ge0$.
\item Hence,  by multiple application of Proposition \ref{Tapia1} (Proposition E.5.2 in \cite{Ta18}), the branch $\tilde{\mu}_j^p(g+tk)$ for any $p \in\{1,2,\dots,\mathfrak{m}\}$ is infinitely many times Gateux differentiable with respect to $g$ and $t$ and the Gateux derivatives are continous. The continuity with respect to $g$ is meant in terms of the $C^1$-uniform topology.
\item Hence,  by multiple application of Proposition \ref{Tapia2} (Proposition E.5.3 in \cite{Ta18}), $\tilde{\mu}_j^p(g+tk)$ is infinitely many times Fr\'{e}chet differentiable with respect to $g$ and $t$ and the Fr\'{e}chet derivatives are continous. The continuity with respect to $g$ is meant in terms of the $C^1$-uniform topology.
\end{itemize}
\end{proof}

\begin{theorem}\label{cor39}
Let $M$ be a compact $m$-dimensional spin manifold which supports a sequence of Riemannian metrics $(h_n)_{n\ge 0}$ such that
\begin{itemize}
\item[$(i)$] the sequence is contained in a compact set of $\Gamma(S^2 T^*M)$ with respect to the uniform $C^1$-topology;
\item[$(ii)$] there exists an $j\geq 0$ such that the eigenvalue $\lambda_{j}(h_n)\rightarrow 0$
for $n\rightarrow +\infty$.
\end{itemize}
Then there exists a Riemannian metric $h$ on $M$ with non-trivial harmonic spinors.
\end{theorem}

\begin{rem}
The metric $h$ is {\it not} the limit of the sequence $(h_n)$, rather it is obtained as an appropriate linear variation of one of its elements.
\end{rem}
\begin{proof}
{\bf Step I.}
\vskip0.2cm\par\noindent
We first need some preliminary observations. Let $g$ be a Riemannian metric on $M$  and $g^t=g+tk$ be the linear variation  determined by a tensor of the form

\begin{equation}
\label{eq:tensor}
k=\|Q\|_{L^{\infty}}^{-1}\;Q,
\end{equation}

where $Q$ is a {\it non-trivial}, i.e., not identically vanishing, symmetric $2$-tensor (possibly depending on $g$) and the $L^{2}$-norm is computed using $g$.

Using a local $g$-orthonormal frame $(e_i)$ of $TM$ which diagonalizes $Q$, it is easy to see that $g^t$ is positive-definite, hence a metric, for all  $|t|<1$.

We will be interested in the case where
$$
Q=Q_{\varphi}
$$
is as in \eqref{eq:Qvar} for some eigenspinor. First note that if $Q_{\varphi}\equiv 0$ for some metric $g$ and a eigenspinor $\varphi\in\Gamma(\Sigma)$ for the corresponding Dirac operator then
\begin{align*}
\left\langle D\varphi,\varphi\right\rangle&=\sum_{i=1}^m\left\langle e_i\cdot\nabla_{e_i}\varphi,\varphi\right\rangle\\
&=\sum_{i=1}^m Q_{\varphi}(e_i,e_i)=0\;,
\end{align*}
so that $\varphi$ would be an harmonic spinor and our result would immediately follow.
We will therefore assume from now on that $Q_{\varphi}$ is non-trivial, {\it for any metric} $g$
{\it and eigenspinor} $\varphi$. In particular, the associated tensors \eqref{eq:tensor} are always well-defined and not identically vanishing. This concludes our preliminary facts.
\vskip0.5cm\par
{\bf Step II.}
\vskip0.2cm\par\noindent
Let $g$ be a metric on $M$ and $g^t=g+tk$ be the linear variation associated to $Q=Q_{\varphi_j^0}$ for some $j\ge0$, with the discrete spectral resolution
$(\varphi_j^t,\lambda_j^t)_{j\ge0}$. As Theorem \ref{BGthm} let $\tilde{\mu}_j^p(g^t)$ be one branch of the eigenvalue $\lambda_j^t$, and the real valued function
\begin{equation}
\label{eq:F}
\begin{split}
F(g,t)&= -2\mathrm{d} \tilde{\mu}_j^p(g^t)/\mathrm{d} t\\
&=\int_M\,g^{t}\left( Q_{\varphi_j^t}, k\right)\mathrm{vol}_{g^t}\\
&=\int_M\,\left<\sum_{i=1}^m e_i\bullet\nabla^t_{(I+tK_g)^{-2}K_g(e_i)}\varphi^t_j,\varphi^t_j\right>^t \text{vol}_{g^t}
\end{split}
\end{equation}
is the derivative of an analytic function, hence analytic in $t$. For the ease of notation we suppress the $j$ and $p$ parameter dependence of the function $F$. Furthermore,
for our special choices of $k$ and $Q$, the function $F$ is locally Lipschitz continuous, hence locally uniformly continuous, in $(g,t)$. The proof of this technical fact is a consequence of following considerations:
\begin{itemize}

\item By Theorem \ref{BGthm}, the branches  of $\lambda_j(g+tk)$ are infinitely many times Fr\'{e}chet derivable with respect to $g$ and $t$ and the Fr\'{e}chet derivatives are continous.  The continuity with respect to $g$ is meant in terms of the $C^1$ uniform topology.
\item We choose now  $k=k_g:=\|Q_{\varphi^0_j}\|_{L^{\infty}}^{-1}\;Q_{\varphi^0_j}$ and compute the first Gateaux derivative of $F$ with respect to $g$:
    \begin{equation}\label{der_F}
    \begin{split}
    &F(g,t) = -2d\tilde{\mu}_j^p(g+tk_g).k_g\\
    &dF(g,t).w = -2d[d\tilde{\mu}_j^p(g+tk_g).k_g].w = \\
    &=-2d^2\tilde{\mu}_j^p(g+tk_g).(d(g+tk_g).w,k_g)-2d\tilde{\mu}_j^p(g+tk_g).(dk_g.w),
    \end{split}
    \end{equation}
where $df.w$ denotes the first Gateaux derivative of $f$ (linear in $w$) and $d^2f.(w,v)$ the second Gateaux derivative of $f$ (bilinear, i.e. linear in both $w$ and $v$). The first  Gateaux derivative of $k_g$ is

\begin{equation}\label{der_k_g}
dk_g= \frac{dQ_{\varphi_j^0}}{\|Q_{\varphi^0_j}\|^2_{L^{2}}}
-\frac{d\|Q_{\varphi_j^0}\|_{L^{\infty}}}{\|Q_{\varphi^0_j}\|^2_{L^{\infty}}}Q_{\varphi_j^0}
\end{equation}

Inserting (\ref{der_k_g}) into (\ref{der_F}) we conclude that the first Gateaux derivative of $F$ with respect to $g$ exists and is continuous in $g$.
\item Hence, the first Gateaux derivative of $F$ with respect to $(g,t)$ exists and is continuous in $(g,t)$. The continuity with respect to $g$ is meant in terms of the $C^1$ uniform topology.
\item Hence, by Proposition \ref{Tapia1} (Proposition E.5.2 in \cite{Ta18}), the first Fr\'{e}chet derivative of $F$ with respect to $(g,t)$ exists, and it is continuous in $(g,t)$.
\item Hence, by proposition \ref{Tapia3} (Proposition E.4.6 in \cite{Ta18}), the function $F=F(g,t)$ is locally Lipschitz continous in $(g,t)$.
\end{itemize}
\vskip0.5cm\par
{\bf Step III.}
\vskip0.2cm\par\noindent
By $(i)$, up to taking a subsequence, we may assume that $(h_n)$ is convergent in the uniform $C^1$-topology; we stress that the limit $h_\infty$ is  not a metric in general but some possibly degenerate symmetric $2$-tensor. Let now
$$\mathcal G=\overline{(h_n)}=(h_n)\cup h_\infty$$
be the closure of $(h_n)$ and define

%

\begin{equation}
\label{eq:restricted}
\mathcal D=\left\{(h_n,t)\left|\,|t|\le \frac{1}{2}\right.\right\}
\end{equation}

\noindent and remark that the function \eqref{eq:F} is well-defined on this set by construction. By $(i)$, the set $\mathcal G$ is compact, hence the closure
\begin{equation}
\label{eq:identitytriangle}
\overline{\mathcal D}=\mathcal D\cup \left\{(h_\infty,t)\left|\, |t|\le \frac{1}{2}\right.\right\}
\end{equation}
of $\mathcal D$ is compact too. We note that identity \eqref{eq:identitytriangle} is a simple consequence of the rectangular form of $\mathcal{D}$. The function \eqref{eq:F} is (locally and hence globally on $\mathcal{D}$) uniformly continuous by Step II, hence it admits a (unique) continuous extension to $\overline{\mathcal D}$, which we still denote by the same symbol $F$. Clearly $F$ is uniformly continuous also on $\overline{\mathcal D}$.
\vskip0.5cm\par
{\bf Step IV.}
\vskip0.2cm\par\noindent
From now on, for simplicity of exposition, we will refer to any element of $\mathcal G$ different from $h_\infty$ directly as a {\it metric $g\in\mathcal G$}. Now

\begin{equation*}
F(g,0)=\frac{\|Q_{\varphi_{j}}\|^2_{L^2}}{\|Q_{\varphi_{j}}\|_{L^{\infty}}}\ge\inf_{g\in\mathcal{G}}\frac{\|Q_{\varphi_{j}}\|^2_{L^2}}{\|Q_{\varphi_{j}}\|_{L^{\infty}}}=:C\ge0
\end{equation*}
for all metrics $g\in\mathcal G$.
%

If $C=0$, since $\mathcal{G}$ is compact and $F(\cdot,0)$ is continuous by Step II, there exists a $g\in\mathcal{G}$ such that $Q_{\varphi_{j}}=0$, the spinor $\varphi_{j}$ is harmonic and we are finished. From now on, we assume that $C>0$.\par
By uniform continuity on $\overline{\mathcal D}$

\begin{equation*}
\label{boundF}
F(g,t)\ge \frac{C}{2},
\end{equation*}
for all $|t|<\delta\le\frac{1}{2}$ for some $\delta$ {\it independent} of $g\in \mathcal G$.


\noindent For this range of $t$, we shall now consider the Taylor expansion at $t=0$ of the $p$-th branch of $j$-th eigenvalue.
For any metric $g\in\mathcal G$, we have
\begin{equation}\label{Taylor}
\tilde{\mu}_j^p(g^t)=\tilde{\mu}_j^p(g) + \alpha(g,\theta)t
\end{equation}
with
\begin{equation*}
\alpha(g,t)=\frac{\mathrm{d} \tilde{\mu}_j^p(g^t)}{\mathrm{d} t}=-\frac{1}{2} F(g,t)\;,
\end{equation*}
where $\theta:(-\delta,\delta)\rightarrow (-\delta,\delta)$ is a continuous function such that $\theta(0)=0$.

We then have
\begin{equation}
\label{eq:range}
\alpha(g,t)\le-\frac{C}{4},
\end{equation}
for any metric $g\in\mathcal G$ and $|t|<\delta$.

\vskip0.2cm\par
Now, observe that $\lambda_j(g)=\tilde{\mu}_j^p(g)=0$ for some metric $g\in\mathcal G$ means that $g$ has non-trivial harmonic spinors. From now on, we therefore assume $\tilde{\mu}_j^p(g)\neq 0$ for all metrics $g\in\mathcal G$.

Given a metric $g\in\mathcal G$ and the associated linear variation $g^t$ for $|t|<\delta$, we need to consider two separate cases:
\begin{itemize}
\item Case $\tilde{\mu}_j^p(g)\ge0$.
\vskip0.1cm\par\noindent

We have
\begin{equation}
\label{dsajkh}
0\le\tilde{\mu}_j^p(g^t)\le \tilde{\mu}_j^p(g)-\frac{C}{4}t,
\end{equation}
for all $0\le t<\delta$. So, the right hand side vanishes if
\begin{equation}
t = \frac{4\tilde{\mu}_j^p(g)}{C}.
\end{equation}
By (ii) we can choose $g\in\mathcal{G}$ such that $t<\delta$ and, therefore, $\tilde{\mu}_j^p(g^t)=0$.

\item Case $\tilde{\mu}_j^p(g)\le0$.
\vskip0.1cm\par\noindent
We have
\begin{equation}
\label{dsajkhII}
\tilde{\mu}_j^p(g)-\frac{C}{4} t\le\tilde{\mu}_j^p(g^t)\le0,
\end{equation}
for all $-\delta< t\le 0$.So, the left hand side vanishes if
\begin{equation}
t = \frac{4\tilde{\mu}_j^p(g)}{C}.
\end{equation}
By (ii) we can choose $g\in\mathcal{G}$ such that $t<\delta$ and, therefore, $\tilde{\mu}_j^p(g^t)=0$.

\end{itemize}
Hence, by Corollary \ref{ev_diff}, the Dirac operator $D^{g^t}$ has an harmonic spinor and the proof is completed.
\end{proof}

\begin{rem} 
It is straightforward to see that Theorem \ref{BGthm}, Corollary \ref{ev_diff} and Theorem \ref{cor39} hold true for a generic Dirac bundle. 
\end{rem}


During the proof of the main Theorem \ref{thm11} we
will need the following auxiliary result, whose proof is straightforward and therefore omitted.
\begin{lemma}\label{cor310}
Let $(M,g)$ be a Riemannian spin manifold and
$g^t=t\,g$
be the linear variation of Riemannian metrics, for all $t>0$. Then
$$
\lambda_j^t=t^{-1/2}\lambda_j\;,
$$
for all $j\ge 0$.
\end{lemma}
\subsection{Berger Metrics}
The one-parameter family of metrics $(g_s)_{s>0}$  on $S^{2k+1}$ known as the {\it Berger metrics} can be conveniently described in terms of the Hopf fibration $S^{2k+1}\rightarrow\mathbf{C}P^k$, $k\ge 1$, where $S^{2k+1}$ is
equipped with the round metric of constant curvature $1$
and $\mathbf{C}P^k$ with the Fubini-Study metric. Recall that the Hopf fibration is a Riemannian submersion with typical fiber $S^1$.
The Berger metric $g_s$ is obtained  rescaling the length of the $S^1$ fibers  by a positive constant $s$ while keeping the metric on the orthogonal complement to the fibers unchanged.
\begin{proposition}(see \cite[pp.~8--9]{Ba96})\label{Berger}
The Berger sphere $(S^{2k+1},g_s)$, $k$ odd,  admits non-trivial harmonic spinors for $s=2(k+1)$.
\end{proposition}

\section{Proof of the Main Theorem}
We split the proof into two steps. First we give conditions under which
a sequence of Riemannian metrics on a compact manifold has a subsequence convergent in the uniform $C^1$-topology. Next, we prove the main theorem about
the existence of harmonic spinors.

Let $(M,g)$ be a compact Riemannian manifold and $\nabla$ its Levi-Civita connection.
For any non-negative integer $k$, we consider the Sobolev space $H^k(M,S^2 T^*M)$ of symmetric $2$-tensors on $M$, with norm given by
\begin{equation}
\label{Sob}
\|h\|^2_{H^k}=\sum_{j=0}^k\int_M|\underbrace{\nabla\nabla\ldots\nabla}_{j-\text{times}}h|^2\mathrm{vol}_{g}
\end{equation}
for all $h\in\Gamma(S^2 T^*M)$. Sobolev norms associated to different Riemannian metrics are all equivalent and induce the same topology of a Hilbert space.
\begin{lemma}\label{limit}
Let $(M,g)$ be a compact Riemannian manifold and $(h_n)$ a sequence of Riemannian metrics bounded in Sobolev norm for a sufficiently large
$k\geq 0$. Then, up to taking a subsequence, $(h_n)$ is contained in a compact set with respect to the uniform $C^1$-topology.
\end{lemma}
\begin{proof}
By Rellich Lemma, up to taking a subsequence, we may assume that $(h_n)$ is convergent in the uniform $C^1$-topology to some $h_{\infty}$. The closure
$$\mathcal G=\overline{(h_n)}=(h_n)\cup\{ h_{\infty}\}$$ is sequentially compact.
\end{proof}
We recall that an oriented hypersurface of a spin manifold inherits a natural Dirac bundle structure (see Proposition \ref{DiracBoundary}).
\begin{theorem}\label{harmonic}
Let $(M,g)$ be a closed connected Riemannian $m$-dimensional spin manifold with spinor bundle $\Sigma\to M$ and $N$ a $0$-codimensional submanifold with an oriented boundary. Assume there exists a Riemannian metric on $\partial N$ with non-trivial harmonic sections of the Dirac bundle $\Sigma|_{\partial N}\to \partial N$. Then there is a Riemannian metric on $M$ with harmonic spinors.
\end{theorem}
\begin{proof}
For any $t>0$ the manifold $M$ is diffeomorphic to the iterated connected sum along a hypersurface $\partial N$:
$$
M^t=(M\setminus \operatorname{Int}(N))\,\#\, U\,\#\, V^t\, \#\, W^t\,\#\,N\;,
$$
where
\begin{equation}
\begin{split}
U&:=[-1,0]\times\partial N ,\\
V^t&:= [0,t]\times\partial N,\\
W^t&:=[t, t+1]\times\partial N .
\end{split}
\end{equation}
In other words, $M$ admits a decomposition as in Proposition \ref{thm27bis}.
We equip $M^t$ with the Riemannian metric $g^t$, defined by
\begin{equation}\label{def_g_t}
g^t=
\begin{cases}
g & \text{on $M\setminus \operatorname{Int}(N)$},\\
\mathrm{d}u^2+\Big(1-\psi(u)+\psi(u)\rho^2(u)\Big)\h &\text{on $U$},\\
\mathrm{d}u^2+\rho^2(u)\h&\text{on $V^t$},\\
\mathrm{d}u^2+\Big(\rho^2(u)-\chi(u)\rho^2(u)+\chi(u)K\Big)\h &\text{on $W^t$},\\
K\,g & \text{on $N$},
\end{cases}
\end{equation}
where $K$ is the constant $\rho^2(t+1)$ for some smooth function $\rho$ to be specified later,
and
$\psi$, $\chi:\mathbf{R}\rightarrow[0,1]$  smooth functions such that
\begin{equation}
\psi(u)=
\begin{cases}
0 & \text{for $u\leq -\frac{3}{4}$}\\
1&\text{for $u\geq -\frac{1}{4}$}
\end{cases}\;,\qquad
\chi(u)=
\begin{cases}
0 & \text{for $u\leq t+\frac{1}{4}$}\\
1&\text{for $u\geq t+\frac{3}{4}$}
\end{cases}\;.
\end{equation}
In order for the metrics on $M\setminus\mathrm{Int}(N)$ and $U$ to join smoothly, it is sufficient to identify a tubular neighbourhood of $\partial (M\setminus\mathrm{Int}(N)) =\partial N$ with an open subset of $U$, following the flow of the outer normal to $\partial N$ as in, e.g., Theorem 9.20 in \cite{Le18}. We note that $\partial N$ is a closed embedded submanifold of $M$.
A similar remark applies to $W^t$ and $N$.

The spin bundle over $(M,g)$ can be stretched over $(M^t,g^t)$.
By Propositions \ref{thm27bis} and \ref{prop55}, any eigenvalue $\lambda_j^2(t)$ of the Dirac Laplacian on $(M^{t},g^{t})$ is dominated by the corresponding Dirichlet eigenvalue $\mu_j(t)$ on  $(V^t,g^{t})$ and there is $_j$ so that
\begin{equation}
\lambda_j^2(t)
\le \mu_j(t)
=\frac{\pi^2}{t^2}\;.
\end{equation}
\par We are going to consider the one-parameter family of metrics
$$\widetilde g^t=\frac{1}{t^{\alpha}}g^t\;,$$
where $0<\alpha<2$. For them the inequality becomes $\widetilde\lambda_j^2(t)\leq\widetilde\mu_j(t)=\pi^2t^{\alpha-2}$ hence
$$
0\leq \widetilde\lambda_j^2(t)\xrightarrow[t\to\infty]{} 0\;.
$$
This implies condition (ii) of Theorem \ref{cor39}.
\par We now specify
$$\rho(u)=\exp\left (-\frac{u}{2(m-1)}\right)$$
and compute the Sobolev norm of the metric $\widetilde g^{t}$ over $V^{t^\alpha}$. Taking the metric $g^1$ as reference, this
Sobolev norm  satisfies for all $k\ge0$ the growth condition
\begin{equation}
\|\tilde{g}^{t^\alpha}\|^2_{\,H^k(V^{t^\alpha})}=\sum_{\left|\alpha\right|\le
k}\int_{V^{t^\alpha}}\,d\text{vol}_{g^1}\left<
{\nabla}^{\alpha}\tilde{g}^{t^\alpha},{\nabla}^{\alpha}\tilde{g}^{t^\alpha}\right>\le\text{const},
\end{equation}
for some positive constant independent of the parameters $t$ and $\alpha$, as long as $\alpha>0$.
The diffeomorphism
\begin{equation}
\Phi^{t^\alpha}:V^1\rightarrow V^{t^\alpha},\quad (u,y)\mapsto(t^\alpha u,y)
\end{equation}
\noindent induces the metric $(\Phi^{t^\alpha})^*\tilde{g}^{t^\alpha}$ on $V^1$ for which,
taking $(\Phi^1)^*g^1$ as reference metric in equation (\ref{Sob})
for the computation of the Sobolev norm,
\begin{equation}
\|(\Phi^{t^\alpha})^*\tilde{g}^{t^\alpha}\|^2_{\,H^k(V^1)}\le\text{const}
\end{equation}
holds true for all $k\ge0$.\\

If $\Psi$ denotes a fixed diffeomorphism mapping $M$ to
$M^1$, then,  by Lemma \ref{cor310} the metrics
$h_j:=(\Psi)^*(\Phi^{t^\alpha})^*\tilde{g}^{t_j^\alpha}$ for a fixed $\alpha\in]0,2[$
satisfy the assumptions of Lemma \ref{limit} and Theorem
\ref{cor39} for a sequence $\{t_j\}_{j\ge 0}$ such that
$t_j\rightarrow +\infty$ as $j\rightarrow\infty$. Hence, there
exists a Riemannian metric on $M$ such that the corresponding Dirac
operator has zero as eigenvalue and the proof is completed.
\end{proof}

\noindent We can proceed now with the proof of the main theorem.
\begin{proof}[Proof of Theorem \ref{thm11}]
Let $p\in M$ be fixed. If $\delta(p)>0$ is the injectivity radius of
$p$, then
\begin{equation}N:=\exp_p\left(\frac{1}{\delta(p)}B^{\mathbf{R}^m}(0,1)\right)\subset M\end{equation} is a geodesic
ball centered at $p$, whose boundary $\partial N$ is diffeomeorphic
to $S^{m-1}$. We analyze several cases:
\begin{itemize}
\item $m=2$: the circle $S^1$ has two spin structures and only one
admits harmonic spinors. Therefore Theorem \ref{harmonic} cannot be
applied. This is in accordance with the fact that $S^2$ has no
harmonic spinors for all Riemannian metrics.
\item $m=3$: as we have just seen the classic Dirac operator on $S^2$
has no vanishing eigenvalues and, again, Theorem \ref{harmonic} cannot
be applied.
\item $m\ge 4$: recalling that there is a unique spin stricture on $S^{m-1}$, we have to distinguish several subcases:
\begin{itemize}
\item $m=4k+4$, for $k\in\mathbf{N}_0$. By Proposition \ref{Berger},
the Berger metric on $S^{m-1}$ admits a non trivial harmonic spinor,
so that Theorem \ref{harmonic} applies and the statement of Theorem
\ref{thm11} holds.
\item  $m=4k+5$, for $k\in\mathbf{N}_0$. We can now apply Theorem \ref{thm11}
to $S^{4k+4}$, as we have just seen, to conclude that there is a
metric on $S^{4k+4}$ admitting non trivial harmonic spinors. We
continue with Theorem \ref{harmonic} and the statement of the Theorem
is proved.
\item $m=4k+6$, for $k\in\mathbf{N}_0$. We apply Theorem \ref{thm11}
to $S^{4k+5}$ and the rest follows analogously.
\item $m=4k+7$, for $k\in\mathbf{N}_0$. We apply Theorem \ref{thm11}
to $S^{4k+6}$ and the rest follows analogously.
\end{itemize}
\end{itemize}
The proof is completed.\\
\end{proof}

\begin{rem} Why does the proof of Theorem \ref{thm11} not work for the Laplace-Beltrami operator or for the complex Laplacian as well?
Of course Theorem \ref{harmonic} holds for general Dirac bundles and
since the constant functions on $S^{m-1}$ are harmonic differential
forms of degree $0$, we can trivially conclude that there are always
non trivial harmonic forms on any Riemannian manifolds, namely the
constant functions. The question now is if it is possible to say
something for differential forms of fixed degree $k>0$.
Unfortunately, both Euler and Clifford operators do not preserve the
graduation of the differential form bundles and there is no version
of Theorem \ref{BGthm} for the Laplace-Beltrami operator or for the
complex Laplacian acting on forms of fixed degree. So, the
constructions in the proof of Theorem \ref{harmonic} cannot be
mimicked for pure forms. If this were the case, we could conclude
for instance that the De Rham cohomology $H_k(S^m,\mathbf{R})$ does
not vanish for $k\neq0,m$, which is of course a contradiction.
 \end{rem}

\section{Conclusion}
We have studied the spectrum of a general Dirac Operator under
variation of the Riemannian metric and extended a result of
Bourguignon and Gauduchon on the first derivative of the eigenvalues
from the special case of the spinor bundle to the general Dirac
bundle case. In conjunction with the extension of the
Courant-Hilbert Theorem for upper bounds of the Dirac Laplacian, we
proved a result on the existence of Riemannian metrics allowing for
harmonic Dirac bundle sections. In particular, we could show that
any closed spin manifold of dimension $m\ge4$ can be always be
provided with a Riemannian metric admitting harmonic spinors. Since
in dimension $m=1,2$, there are counterexamples, the conjecture
remains open for $m=3$.

\appendix
\section{Some Results from Functional Analysis}\label{app}

\begin{theorem}[Rellich]\label{Rellich}
Let $T(x)$ be a selfadjoint holomorphic family of type (A) defined for $x$ in a neighbourhood of an interval $I_0$ of the real axis. Furthermore, let $T(x)$ have compact resolvent. Then, all eigenvalues of $T(x)$ can be represented by functions which are holomorphic on $I_0$. More precisely, there is a sequence of scalar-valued functions $\mu_n(x)$ and a sequence of vector-valued functions $\varphi_n(x)$, all holomorphic on $I_0$, such that if $x\in I_0$, the $\mu_n(x)$ represent all the repeated eigenvalues of $T(x)$ and the $\varphi_n(x)$ form a complete orthonormal family of the associated eigenvectors of $T(x)$.
\end{theorem}

\begin{proposition}\label{Tapia1}
Let $X,Y$ be normed linear spaces. Assume that $f:X\rightarrow Y$ has a directional variation in $X$, i.e. there exists a directional variation $f^{\prime}(x)(\eta)$ for all $x,\eta\in X$. Assume further that:
\begin{itemize}
\item[(i)] For fixed $x$, $f^{\prime}(x)(\eta)$ is continuous in $\eta$ at $\eta=0$.
\item[(ii)] For fixed $\eta$, $f^{\prime}(x)(\eta)$ is continuous in $x$ for all $x\in X$.
\end{itemize}
Then, $f^{\prime}(x)\in\mathcal{L}(X,Y)$, i.e. $f^{\prime}(x)$ is a bounded linear operator.
\end{proposition}

\begin{proposition}\label{Tapia2}
Consider $f:X\rightarrow Y$, where $X$ and $Y$ are normed linear spaces. Suppose that $f$ is Gateaux
differentiable in an open set $U\subset X$. If $f^{\prime}(x)$ is continuous at $x\in U$, then $f^{\prime}(x)$  is a Fr\'{e}chet derivative.
\end{proposition}
\begin{proposition}\label{Tapia3}
Consider $f:X\rightarrow Y$, where $X$ and $Y$ are normed linear spaces.
If $f$ is Fr\'{e}chet differentiable at $x\in X$, then f is continuous at $x$. If the Fr\'{e}chet derivative $f^{\prime}(x)$ is continuous for all $x\in X$, then $f$ is locally Lipschitz continuous.
\end{proposition}

\section*{Acknowledgements}
I would like to express my gratitude to Guido Franchetti and Andrea Santi for their very valuable contribution to some very technical parts of this paper and its predecessors, and to Bernd Ammann and Urs Semmelmann for their precious help to clarify several questions around the differentiability of the Dirac eigenvalues. All possibly remaining mistakes are mine.

\end{document}